\documentclass[dvips,preprint, 11pt]{imsart}
\usepackage{amsmath}
\usepackage{amsfonts}
\usepackage{mathrsfs}

\usepackage{ amsmath, amsfonts, amssymb, amsthm, amscd}
\usepackage[T1]{fontenc}
\usepackage[english]{babel}
\usepackage{color}
\usepackage{pgf,pgfarrows,pgfnodes,pgfautomata,pgfheaps}
\usepackage{amsmath,amssymb}
\usepackage[latin1]{inputenc}

\usepackage[latin1]{inputenc}
\usepackage{graphicx}

\fontsize{11pt}{14.5pt}\selectfont

\def\Dy#1{\Frac{\partial #1}{\partial y}}

\def\Dy_1y_1#1{\Frac{\partial^2 #1}{\partial y_1^2}}

\newtheorem{Theorem}{Theorem}[part]
\newtheorem{Definition}{Definition}[part]
\newtheorem{Proposition}{Proposition}[part]
\newtheorem{Assumption}{Assumption}[part]
\newtheorem{Lemma}{Lemma}[part]

\newtheorem{Remark}{Remark}[part]

\makeatletter \@addtoreset{equation}{section}

\@addtoreset{Definition}{section}

\@addtoreset{Theorem}{section}

\@addtoreset{Proposition}{section}

\@addtoreset{Assumption}{section}

\@addtoreset{Corollary}{section}

\@addtoreset{Lemma}{section}

\@addtoreset{Remark}{section}

\@addtoreset{Example}{section}

\makeatother \makeatletter

\def \Int{\displaystyle\int}
\def \Frac{\displaystyle\frac}

\def \Sup{\displaystyle\sup}

\def \be{\begin{eqnarray}}
\def \ee{\end{eqnarray}}
\def \b*{\begin{eqnarray*}}
\def \e*{\end{eqnarray*}}

\def \E{\mathbb{E}}

\def \P{\mathbb{P}}
\def \R{\mathbb{R}}

\def \[{[\,\!\![}
\def \]{]\,\!\!]}

\def \1{{\bf 1}}

\begin{document}

\begin{frontmatter}

\title{ Large Deviation Principles of Obstacle Problems  for Quasilinear Stochastic PDEs}
\runtitle{Quasilinear SPDEs}

\author{\fnms{Anis} \snm{Matoussi}\ead[label=e1]{anis.matoussi@univ-lemans.fr}}
\thankstext{t2}{The first author
 was  partially supported by  chaire Risques Financiers de la fondation du risque, CMAP-\'Ecole Polytechniques, Palaiseau-France  }
\address{University of Maine \\
Risk and Insurance Institute of Le Mans \\
Laboratoire  Manceau de Math\'ematiques\\ Avenue Olivier Messiaen\\
 \printead{e1}}
  \author{\fnms{Wissal} \snm{Sabbagh}\corref{}\ead[label=e2]{wissal.sabbagh@univ-lemans.fr}}
\address{ University of Maine \\
Risk and Insurance Institute of Le Mans \\
Laboratoire Manceau de Math\'ematiques\\ Avenue Olivier Messiaen \\ \printead{e2}}
 \affiliation{University of Le Mans}

\author{\fnms{Tusheng} \snm{Zhang}\corref{}\ead[label=e3]{tusheng.zhang@manchester.ac.uk}}
\address{ School of Mathematics, University of Manchester, \\
Oxford Road, Manchester M13 9PL, England, UK \\
 \printead{e3}}
 \affiliation{University of Manchester}


\runauthor{Matoussi, Sabbagh and Zhang }


\begin{abstract}
In this paper, we present a sufficient condition for  the large deviation criteria of Budhiraja, Dupuis and Maroulas for functionals of  Brownian motions. We then establish a large deviation principle for obstacle problems  of quasi-linear stochastic partial differential equations. It turns out that the backward stochastic differential equations will play an important role.
\end{abstract}

\end{frontmatter}

\section{Introduction}
Consider the following  obstacle problems for quasilinear stochastic partial differential equations (SPDEs) in $\R^d$:
\small
\begin{eqnarray}\label{SPDE1}
dU(t,x)&+&\frac{1}{2}\Delta U(t,x)+\sum_{i=1}^d\partial_ig_i(t,x,U(t,x),\nabla U(t,x))dt +f(t,x,U(t,x),\nabla U(t,x))dt\nonumber\\
 &+ &\sum_{j=1}^{\infty} h_j(t,x,U(t,x),\nabla U(t,x))dB^j_t =-R(dt,dx),\\
 U(t,x)&\geq & L(t,x), \quad \quad (t,x)\in \R^+\times \R^d,\nonumber\\
 U(T,x)&=&\Phi (x), \quad\quad x\in \R^d,
 \end{eqnarray}
where $B^j_t, j=1,2,...$ are independent real-valued standard Brownian motions, the stochastic integral against Brownian motions is interpreted as the backward Ito integral, $\Delta $ is the Laplacian operator, $f, g_i, h_j$ are appropriate measurable functions specified later, $L(t,x)$ is the given barrier function, $R(dt, dx)$ is a random measure which is a part of the solution pair $(U, R)$. The random measure $R$ plays a similar role as a local time which prevents the solution $U(t,x)$ from falling below the barrier $L$.
\vskip 0.3cm
Such SPDEs appear in various applications like pathwise stochastic control problems, the Zakai equations in filtering and stochastic control with partial observations.
Existence and uniqueness of the above stochastic obstacle problems were established in \cite{DMZ} based on an analytical approach. Existence and uniqueness of the  obstacle problems for quasi-linear SPDEs on the whole space $\R^d$ and driven by finite dimensional Brownian motions  were studied in \cite{MS} using  the approach of backward stochastic differential equations (BSDEs). Obstacle problems for nonlinear stochastic heat equations driven by space-time white noise were studied by several authors, see \cite{NP},\cite{XZ} and references therein.
\vskip 0.3cm
 In this paper, we are concerned with the small small noise large deviation principle(LDP) of the following  obstacle problems for quasilinear SPDEs:
\small
\begin{eqnarray}\label{SPDE1-1}
dU^{\varepsilon}(t,x)&+&\frac{1}{2}\Delta U^{\varepsilon}(t,x)+\sum_{i=1}^d\partial_ig_i(t,x,U^{\varepsilon}(t,x),\nabla U^{\varepsilon}(t,x))dt +f(t,x,U^{\varepsilon}(t,x),\nabla U^{\varepsilon}(t,x))dt\nonumber\\
 &+ &\sqrt{\varepsilon}\sum_{j=1}^{\infty} h_j(t,x,U^{\varepsilon}(t,x),\nabla U^{\varepsilon}(t,x))dB^j_t =-R^{\varepsilon}(dt,dx),\\
 U^{\varepsilon}(t,x)&\geq & L(t,x), \quad \quad (t,x)\in \R^+\times \R^d,\nonumber\\
 U^{\varepsilon}(T,x)&=&\Phi (x), \quad\quad x\in \R^d.
 \end{eqnarray}
Large deviations for stochastic evolution equations and stochastic
 partial differential equations driven by Brownian motions have been
investigated in many papers, see e.g.\ \cite{Duan-Millet}, \cite{L}, \cite{Manna-Sritharan-Sundar}, \cite{Rockner-Zhang-Zhang}, \cite{CW}, \cite{CR}, \cite{S}, \cite{Budhiraja-Dupuis-Maroulas-1.}, \cite{C} and references therein.
\vskip 0.3cm
 To obtain the large deviation principle, we will adopt the weak convergence approach introduced by Budhiraja, Dupuis and Maroulas in \cite{Budhiraja-Dupuis-Maroulas.} , \cite{Budhiraja-Dupuis-Maroulas-1.} and \cite{Budhiraja-Dupuis}). This approach is now a powerful tool which has been applied by many  people to prove large deviation principles  for various dynamical systems  driven by Gaussian noises, see e.g. \cite{Budhiraja-Dupuis}, \cite{Duan-Millet}, \cite{L}, \cite{Manna-Sritharan-Sundar}, \cite{Rockner-Zhang-Zhang},\cite{Budhiraja-Dupuis-Maroulas-1.}.
 \vskip 0.3cm
 In order to apply the weak convergence approach to the obstacle problems, we first provide a  simple sufficient condition to verify the criteria of Budhiraja-Dupuis-Maroulas. The sufficient condition is particularly suitable for stochastic dynamics generated by stochastic differential equations and stochastic partial differential equations. The important part of the work is to study the so called skeleton equations which, in our case, are the deterministic obstacle problems driven by the elements in the Cameron-Martin space of the Brownian motions. We need to show that if the driving signals converge weakly in the Cameron-Martin space, then the corresponding solutions of the skeleton equations also converge in the appropriate state space. This turns out to be hard because of the singularity caused by the obstacle. To overcome the difficulties, we have to appeal to the penalized approximation of the skeleton equation and to establish some uniform estimate for  the solutions of the approximating equations with the help of  the backward stochastic differential equation representation of the solutions. This is purely due to the technical reason because primarily the LDP problem has not much to do with backward stochastic differential equations.
\vskip 0.3cm
The rest of the paper is organized as follows. In Section 2, we introduce the stochastic obstacle problem and the precise framework. In Section 3, we recall the weak convergence approach of large deviations and present a sufficient condition. Section 4 is devoted  to the study of skeleton obstacle problems. We will show that the solution of the skeleton problem is continuous with respect to the driving signal. The proof of the large deviation principle is in Section 5.

\section{The framework}
\subsection{Obstacle problems}
\label{section:SPDE}
Let  $H:=\mathbf{L}^2\left(\R^d\right)$ be the Hilbert space of square integrable functions with respect to the Lebesgue measure on $\R^d$. The associated scalar product and the norm  are denoted by
$$ (u,v)=\Int_{\R^d}u\left( x\right) v\left(
x\right)dx, \quad\quad\quad |u|=\left (\Int_{\R^d}u^2(x)dx\right )^{1/2}. $$
Let $V:=H(\R^d)$ denote the first order Sobolev space, endowed with the norm and the inner product:
$$\|u\|=\left (\Int_{\R^d}|\nabla u|^2(x)dx+\Int_{\R^d}|u|^2(x)dx\right )^{1/2},$$
$$ <u, v>=\Int_{\R^d}(\nabla u)\cdot (\nabla v)(x)dx+\Int_{\R^d}u(x)v(x)dx.$$
$V^*$ will denote the dual space of $V$. When causing no confusion, we also use $<u,v>$ to denote the dual pair between $V$ and $V^*$.

Our evolution problem will be considered over a fixed time interval
$[0,T]$. Now we introduce the following assumptions.
\begin{Assumption} \label{assgener}
\begin{itemize}
\item[\rm{(i)}] $f:[0,T]\times \R^{d}\times \R\times \R^{d}\rightarrow\R$, $h=(h_1,...,h_i,...):[0,T]\times \R^{d}\times \R\times \R^{d}\rightarrow\R^{\infty}$ and $g=(g_1,...,g_d):[0,T]\times \R^{d}\times\R \times\R^{d}\rightarrow\R^{d}$  are measurable in $(t,x,y,z)$ and
satisfy $ f^0, h^0, g^0 \in \mathbf{L}^2\left( [0,T] \times \R^{d}\right)\cap \mathbf{L}^{\infty}\left( [0,T] \times \R^{d}\right) $ where $f^0(t,x) := f(t,x,0, 0)$, $h^0(t,x) := (\sum_{j=1}^{\infty}h_j(t,
,x,0, 0)^2)^{\frac{1}{2}}$ and $g^0(t,x) := (\sum_{j=1}^{d}g_j(t,
,x,0, 0)^2)^{\frac{1}{2}}$ .
 \item [\rm{(ii)}] There exist constants $c>0$,\,$0<\alpha<1$ and $0<\beta<1$ such that for any $(t,x)\in[0,T]\times\R^d~;~(y_1,z_1),(y_2,z_2)\in\R\times\R^{d}$
\b*
|f(t,x,y_1,z_1)-f(t,x,y_2,z_2)| &\leq & c\big(|y_1-y_2|+|z_1-z_2|\big)\\
\left (\sum_{i=1}^{\infty} |h_i(t,x,y_1,z_1)-h_i(t,x,y_2,z_2)|^2\right )^{1/2} &\leq & c|y_1-y_2|+\beta |z_1-z_2|\\
\left (\sum_{i=1}^{d} |g_i(t,x,y_1,z_1)-g_i(t,x,y_2,z_2)|^2\right )^{1/2} &\leq & c|y_1-y_2|+\alpha |z_1-z_2|.
\e*
\item [\rm{(iii)}]There exists a function $\bar{h}\in L^2(\R^{d})\cap L^{\infty}(\R^{d})$ such that for $(t,x,y,z) \in [0,T]\times \R^{d}\times \R\times \R^{d}$,
$$\left (\sum_{i=1}^{\infty} |h_i(t,x,y,z)|^2\right )^{1/2}\leq \bar{h}(x).$$
\item [\rm{(iv)}]The contract property: $\alpha+\Frac{\beta^2}{2}<\Frac{1}{2}$.
\item [\rm{(v)}]The barrier function $L(t,x):\R^d\rightarrow \R$ satisfies
$$\frac{\partial L(t,x)}{\partial t}, \quad \nabla L(t,x), \quad \Delta L(t,x)\in L^2([0,T]\times \R^d)\cap L^{\infty}([0,T]\times \R^d),$$
where the gradient $\nabla$ and the Laplacian $\Delta$ act on the space variable $x$.
\end{itemize}
\end{Assumption}
\vskip 0.4cm

Let $H_T:=C([0, T], H)\cap L^2([0, T], V)$ be the Banach space endowed with the norm
$$\|u\|_{H_T}=\sup_{0\leq t\leq T}|u_s|+\left (\int_0^T\|u_s\|^2ds\right )^{1/2}.$$
We denote by  $ {\mathcal H}_T$ the space of  predictable, processes  $(u_t, t\geq 0 ) $ such that $u\in H_T$ and that
$$
 \E \,  \big[\underset{ 0 \leq s \leq T}{\Sup} |u_s |_2^2 +   \Int_0^T  \|u_s\|^2 ds\big]<\infty.
$$
The space of test functions  is
$\mathcal{D}= \mathcal{C}_c^{\infty} \left( \R^+\right) \otimes
\mathcal{C}_c^{\infty} \left( \R^d\right)$, where
$\mathcal{C}_c^{\infty} \left(\R^+\right)$ denotes the space of real-valued infinitely differentiable functions with compact supports in $\R^+$
and  $\mathcal{C}_c^{\infty}\left( \R^d\right)$ is the space of
infinitely differentiable functions with compact supports in
$\R^d$.
\vskip 0.3cm
\noindent\\
We  now precise  the definition  of solutions for the reflected quasilinear SPDE (\ref{SPDE1}):
\begin{Definition}
\label{o-pde}We say that a pair $(U, R )$ is a solution of the obstacle problem (\ref{SPDE1}) if
\begin{itemize}
\item[(1)]$U\in {\mathcal H}_T$, $U(t,x)\geq L(t,x)$, $dP\otimes dt\otimes dx$-a.e. and  $U(T,x)=\Phi(x)$, $dx-a.e.$
\item[(2)] $R$ is a random regular measure on $[0,T)\times \R^d$,
\item[(3)] for every $\varphi \in \mathcal D$%
\begin{align}\label{OPDE}
\nonumber &(U_t,\varphi_t)-(\Phi,\varphi_T)+\int_t^T(U_s, \partial_s\varphi_s)ds+\frac{1}{2}\int_t^T<\nabla U_s, \nabla\varphi_s>ds \\
\nonumber &=\int_t^T(f_s(U_s,\nabla U_s),\varphi_s)ds+\sum_{j=1}^{\infty}\int_t^T (h_s^j(U_s, \nabla U_s), \varphi_s)dB_s^j\ \\
&-\sum_{i=1}^d \int_t^T(g_s^i(U_s, \nabla U_s), \partial_i\varphi_s)ds+\int_t^T\int_{\R^d}\varphi_s(x)R(dx,ds),
\end{align}
\item[(4)] $U$ admits a quasi-continuous version $\tilde{U}$, and
$$\int_0^T\int_{\R^d}(\tilde{U}(s,x)-L(s,x))R(dx,ds)=0\quad\quad a.s.$$
\end{itemize}
\end{Definition}
\begin{Remark}
We refer the reader to \cite{DMZ} for the precise definition of regular measures and quasi-continuity of functions on the space $[0, T]\times \R^d$.
\end{Remark}
Let us recall the following result from \cite{MS} and \cite{DMZ}.
\begin{Theorem}
\label{existence:RSPDE}
Let  Assumption 2.1   hold and assume $\Phi(x)\geq L(T,x)$ $dx$-a.e.. Then
there exists a unique solution $(U, R)$ to the obstacle problem (\ref{SPDE1}).
\end{Theorem}
\subsection{The measures $\P^m$}
The operator $ \partial_t + \frac{1}{2} \Delta $, which represents the main linear part in the equation \eqref{SPDE1}, is associated with the Bownian motion in $\mathbb{R}^d$.   The sample space of the Brownian motion is $ \Omega' = \mathcal{C }\left([0, \infty ); \mathbb{R}^d \right)$, the canonical process $(W_t)_{t \geq 0}$ is defined by $ W_t (\omega) = \omega (t)$, for any $ \omega \in \Omega'$, $t \geq 0$ and  the shift operator, $ \theta_t  \, : \,  \Omega' \longrightarrow  \Omega'$, is defined by $ \theta_t (\omega) (s) = \omega (t+s)$, for any $s \geq 0$ and $ t \geq 0$. The canonical filtration $ \mathcal{F}_t^W = \sigma \left( W_s; s \leq t \right)$ is completed by the standard procedure with respect to the probability measures produced by the transition function
$$
P_t (x, dy) = q_t (x-y) dy, \quad t >0, \quad x \in \mathbb{R}^d,
$$
where $ q_t (x) = \left(2\pi t\right)^{- \frac{d}{2}} \exp \left( - |x|^2/2t \right)$ is the Gaussian density. Thus we get a continuous Hunt process $\left(\Omega', W_t, \theta_t, \mathcal{F}, \mathcal{F}^W_t, \mathbb{P}^x \right)$. We shall also use the backward filtration of the future events $ \mathcal{F}'_t = \sigma \left(W_s; \; \,  s \geq t \right)$ for $t\geq 0$. $\mathbb{P}^0$ is the Wiener measure, which is supported by the set $ \Omega'_0 = \{ \omega \in \Omega', \; \, \omega(0) =0 \}$. We also set $ \Pi_0 (\omega) (t) = \omega (t) - \omega (0),\,  t \geq 0$, which defines a map $ \Pi_0 \, : \, \Omega' \rightarrow \Omega'_0$. Then  $\Pi = (W_0, \Pi_0 ) \, : \,  \Omega' \rightarrow \mathbb{R}^d \times \Omega'_0$ is a bijection. For each probability measure on $\mathbb{R}^d$, the probability $\mathbb{P}^{\mu}$ of the Brownian motion started with the initial distribution $\mu$ is given by $$ \mathbb{P}^{\mu} = \Pi^{-1} \left(\mu \otimes \mathbb{P}^0 \right).$$
In particular, for the Lebesgue measure in $\mathbb{R}^d$, which we denote by $ m = dx$, we have
$$ \mathbb{P}^{m} = \Pi^{-1} \left(dx\otimes \mathbb{P}^0 \right).$$
Notice that $\{W_{t-r}, \mathcal{F}'_{t-r}, r\in [0, t]\}$ is a backward local martingale under $\P^m$.
Let $J(\cdot,\cdot): [0, \infty )\times \R^d \rightarrow \R^d$ be a measurable function such that $J\in \mathbf{L}^2\left(
[0,T] \times \mathbb{{R}}^d\rightarrow \mathbb{{R}}^d\right) $ for every $T>0$. We recall the forward and backward stochastic integral defined in \cite{S}, \cite{MS} under the measure $\P^m$.
$$\Int_s^t J(r,W_{r})\ast dW_r=\Int_s^t <J(r,W_{r}), dW_r>+\Int_s^t <J(r,W_{r}), d\overleftarrow{W}_r>. $$
When $J$ is smooth, one has
\begin{equation}\label{forward-backward}
\Int_s^t J(r,W_{r})\ast dW_r=-2\Int_s^t div(J(r,\cdot))(W_r)dr.
\end{equation}
We refer the reader to \cite{MS}, \cite{S} for more details.
\section{A sufficient condition for LDP}
In this section we will recall the criteria obtained in \cite{Budhiraja-Dupuis} for proving a large deviation principle and we will also present a sufficient condition to verify the criteria.
\vskip0.3cm
Let $\mathcal{E}$ be a Polish space with the Borel $\sigma$-field $\mathcal{B}(\mathcal{E})$. Recall
    \begin{Definition}\label{Dfn-Rate function}
       \emph{\textbf{(Rate function)}} A function $I: \mathcal{E}\rightarrow[0,\infty]$ is called a rate function on
       $\mathcal{E}$,
       if for each $M<\infty$, the level set $\{x\in\mathcal{E}:I(x)\leq M\}$ is a compact subset of $\mathcal{E}$.
         \end{Definition}
    \begin{Definition}
       \emph{\textbf{(Large deviation principle)}} Let $I$ be a rate function on $\mathcal{E}$.  A family
       $\{X^\varepsilon\}$ of $\mathcal{E}$-valued random elements is  said to satisfy a large deviation principle on $\mathcal{E}$
       with rate function $I$ if the following two claims
       hold.
       \begin{itemize}
         \item[$(a)$](Upper bound) For each closed subset $F$ of $\mathcal{E}$,
              $$
                \limsup_{\varepsilon\rightarrow 0}\varepsilon\log\mathbb{P}(X^\varepsilon\in F)\leq- \inf_{x\in F}I(x).
              $$
         \item[$(b)$](Lower bound) For each open subset $G$ of $\mathcal{E}$,
              $$
                \liminf_{\varepsilon\rightarrow 0}\varepsilon\log\mathbb{P}(X^\varepsilon\in G)\geq- \inf_{x\in G}I(x).
              $$
       \end{itemize}
    \end{Definition}

\vskip0.3cm
\subsection{A criteria of Budhiraja-Dupuis}

The Cameron-Martin space associated with the Brownian motion $\{B_t=(B_t^1,...,B^j_t,...), t\in[0,T]\}$ is isomorphic to the Hilbert space $K:=L^2([0,T]; l^2)$ with the inner product:
$$
 \langle h_1, h_2\rangle_{K}:=\int_0^T\langle  h_1(s), h_2(s)\rangle_{l^2}ds,
 $$
where
$$l^2=\{a=(a_1,...,a_j,...); \sum_{i=1}^{\infty}a_i^2<\infty\}.$$
$l^2$ is a  Hilbert space with inner product $\langle  a, b\rangle_{l^2}=\sum_{i=1}^{\infty}a_ib_i$ for $a,b\in l^2$.

Let $\tilde{K}$ denote  the class of $l^2$-valued $\{\mathcal{F}_t\}$-predictable processes $\phi$ that belong to the space $K$ a.s..
Let $S_N=\{k\in K; \int_0^T\|k(s)\|_{l^2}^2ds\le N\}$. The set $S_N$ endowed with the weak topology is a compact Polish space.
Set $\tilde{S}_N=\{\phi\in \tilde{K};\phi(\omega)\in S_N, \mathbb{P}\text{-a.s.}\}$.

\vskip0.3cm

 The following result was proved in \cite{Budhiraja-Dupuis}.

\begin{Theorem}\label{thm BD}{ For $\varepsilon>0$, let $\Gamma^\varepsilon$ be a measurable mapping from $C([0,T];\mathbb{R}^\infty)$ into $\mathcal{E}$.
Set $X^\varepsilon:=\Gamma^\varepsilon(B(\cdot))$. Suppose that
there  exists a measurable map $\Gamma^0:C([0,T];\mathbb{R}^\infty)\rightarrow \mathcal{E}$ such that
\begin{itemize}
   \item[(a)] for every $N<+\infty$ and any family $\{k^\varepsilon;\varepsilon>0\}\subset \tilde{S}_N$ satisfying that $k^\varepsilon$ converges in law as $S_N$-valued random elements to some element $k$ as $\varepsilon\rightarrow 0$,
    $\Gamma^\varepsilon\left(B(\cdot)+\frac{1}{\sqrt\varepsilon}\int_0^{\cdot}k^\varepsilon(s)ds\right)$ converges in law to $\Gamma^0(\int_0^{\cdot}k(s)ds)$ as $\varepsilon\rightarrow 0$;
   \item[(b)] for every $N<+\infty$, the set
   $$
 \left\{\Gamma^0\left(\int_0^{\cdot}k(s)ds\right); k\in S_N\right\}
  $$
   is a compact subset of $\mathcal{E}$.
 \end{itemize}
Then the family $\{X^\varepsilon\}_{\varepsilon>0}$ satisfies a large deviation principle in $\mathcal{E}$ with the rate function $I$ given by
\begin{equation}\label{rate function}
I(g):=\inf_{\{k\in K;g=\Gamma^0(\int_0^{\cdot}k(s)ds)\}}\left\{\frac12\int_0^T\|k(s)\|_{l^2}^2ds\right\},\ g\in\mathcal{E},
\end{equation}
with the convention $\inf\{\emptyset\}=\infty$.
 }
 \end{Theorem}
 \subsection{A sufficient condition}
 Here is a sufficient condition for verifying the assumptions in Theorem 3.1 .
 \begin{Theorem}For $\varepsilon>0$, let $\Gamma^\varepsilon$ be a measurable mapping from $C([0,T];\mathbb{R}^\infty)$ into $\mathcal{E}$.
Set $X^\varepsilon:=\Gamma^\varepsilon(B(\cdot))$. Suppose that
there  exists a measurable map $\Gamma^0:C([0,T];\mathbb{R}^\infty)\rightarrow \mathcal{E}$ such that
\begin{itemize}
   \item[(i)] for every $N<+\infty$, any family $\{k^\varepsilon;\varepsilon>0\}\subset \tilde{S}_N$ and any $\delta>0$,
    $$\lim_{\varepsilon\rightarrow 0}P(\rho \left (Y^\varepsilon, Z^\varepsilon\right )>\delta )=0,$$
    where $Y^\varepsilon=\Gamma^\varepsilon\left(B(\cdot)+\frac{1}{\sqrt\varepsilon}\int_0^{\cdot}k^\varepsilon(s)ds\right)$,  $Z^\varepsilon=\Gamma^0\left(\int_0^{\cdot}k^\varepsilon(s)ds\right)$ and  $\rho(\cdot, \cdot)$ stands for the metric in the space $\mathcal{E}$
    \item[(ii)] for every $N<+\infty$ and any family $\{k^\varepsilon;\varepsilon>0\}\subset {S}_N$ satisfying that $k^\varepsilon$ converges weakly to some element $k$ as $\varepsilon\rightarrow 0$, $\Gamma^0\left(\int_0^{\cdot}k^\varepsilon(s)ds\right)$ converges  to $\Gamma^0(\int_0^{\cdot}k(s)ds)$ in the space
   $\mathcal{E}$.
 \end{itemize}
Then the family $\{X^\varepsilon\}_{\varepsilon>0}$ satisfies a large deviation principle in $\mathcal{E}$ with the rate function $I$ given by
\begin{equation}\label{rate function}
I(g):=\inf_{\{k\in K;g=\Gamma^0(\int_0^{\cdot}k(s)ds)\}}\left\{\frac12\int_0^T\|k(s)\|_{l^2}^2ds\right\},\ g\in\mathcal{E},
\end{equation}
with the convention $\inf\{\emptyset\}=\infty$.
 \end{Theorem}
 \begin{Remark}. When proving a small noise large deviation principle for stochastic differential equations/stochastic partial differential equations, condition (i) is usually not difficult to check because the small noise disappears when $\varepsilon\rightarrow 0$.
 \end{Remark}

 \vskip 0.4cm
 \noindent{\bf Proof}. We will  show that the conditions in Theorem 3.1 are fulfilled. Condition (b) in Theorem 3.1 follows from condition (ii) because $S_N$ is compact with respect to the weak topology. Condition (i) implies that for any bounded, uniformly continuous function $G(\cdot)$ on $\mathcal{E}$,
 $$\lim_{\varepsilon \rightarrow 0}E[|G(Y^\varepsilon)-G(Z^\varepsilon)|]=0.$$
 Thus, condition (a) will be satisfied if  $Z^\varepsilon $ convergence in law to $\Gamma^0(\int_0^{\cdot}k(s)ds)$ in the space $\mathcal{E}$. This is indeed true since
the mapping $\Gamma^0$ is continuous by condition (ii) and $k^\varepsilon$ converge in law as $S_N$-valued random elements to $k$.
The proof is complete.

\section{Skeleton equations}
Recall  $K:=L^2([0, T], l^2)$.  Let $k\in K$ and consider the deterministic obstacle problem:
\begin{eqnarray}\label{2.1}
du (t,x)&+&\frac{1}{2}\Delta u(t,x)+\sum_{i=1}^d\partial_ig_i(t,x, u(t,x), \nabla u(t,x))dt +f(t,x,u(t,x), \nabla u(t,x))dt\nonumber\\
 &+ &\sum_{j=1}^{\infty} h_j(t,x,u(t,x), \nabla u(t,x))k^j_tdt =-\nu(dt,dx),\\
 u(t,x)&\geq & L(t,x), \quad \quad (t,x)\in \R^+\times \R^d,\nonumber\\
 u(T,x)&=&\Phi (x), \quad\quad x\in \R^d.
 \end{eqnarray}
 Denote by $u^{\varepsilon}$ the solution of equation (\ref{2.1}) with $k^{\varepsilon}$ replacing $k$. The main purpose of this section is to show that  $u^{\varepsilon}$ converges to $u$ in the space $H_T$ if $k^{\varepsilon}\rightarrow k$ weakly in the Hilbert space $K$.
To this end, we need to establish a number of preliminary results.
\vskip 0.3cm
Consider the penalized equation:
\begin{eqnarray}\label{2.2}
du^n (t,x)&+&\frac{1}{2}\Delta u^n(t,x)+\sum_{i=1}^d\partial_ig_i(t,x, u^n(t,x), \nabla u^n(t,x))dt +f(t,x,u^n(t,x), \nabla u^n(t,x))dt\nonumber\\
 &+ &\sum_{j=1}^{\infty} h_j(t,x,u^n(t,x), \nabla u^n(t,x))k^j_tdt =-n(u^n(t,x)-v(t,x))^-dt ,\\
 u^n(T,x)&=&\Phi (x), \quad\quad x\in \R^d.
 \end{eqnarray}
 \vskip 0.4cm
 For later use, we need to show that for any $M>0$,  $u^n\rightarrow u$ uniformly over the bounded subset $\{k; \|k\|_K\leq M\}$ as $n\rightarrow \infty$. For this purpose, it turns out that we have to appeal to the BSDE representation of the solutions.
  Let $Y^n_t:=u^n(t, W_t)$, $Z^n_t=\nabla u^n(t, W_t)$. Then it was shown in \cite{MS} that
 $(Y^n, Z^n)$ is the solution of the backward stochastic differential equation under $\P^m$:
 \begin{eqnarray}\label{2.4}
Y_{t}^{n}&=&\Phi(W_{T})+\Int_{t}^{T}f(r,W_{r},Y_{r}^{n},Z_{r}^{n})dr+\sum_{j=1}^{\infty}\Int_{t}^{T}h_j(r,W_{r},Y_{r}^{n},Z_{r}^{n})k_r^jdr\nonumber\\
&&+n\Int_t^T (Y^n_r-S_r)^-dr +\frac{1}{2}\Int_s^T g(r,W_{r},Y_{r}^{n},Z_{r}^{n})\ast dW_r-\Int_{t}^{T}Z_{r}^{n}dW_{r}.\nonumber\\
& &
\end{eqnarray}
Where $S_r=L(r, W_r)$ satisfies
\begin{equation}\label{2.5}
 dS_r=\frac{\partial L}{\partial r}(r,W_r)dr+ \frac{1}{2}\Delta L(r,W_r)dr+\nabla L(r,W_r)dW_{r}.
\end{equation}
The following result is a uniform estimate for $(Y^n, Z^n)$.
\begin{Lemma}
For $M>0$, we have the following estimate:
\begin{eqnarray}\label{2.6}
&& \underset{\{k\in K; \|k\|_K\leq M\}}{\Sup}\underset{n}{\Sup}\left \{ \E^m[\underset{0\leq t\leq T}{\Sup}|Y_{t}^{n}|^2]+\E^m[\Int_0^T|Z_t^n|^2dt]+\E^m[\left (n\Int_0^T(Y_t^n-S_t)^-dt\right )^2]\right \}\nonumber\\
&\leq & c_M \left [|\Phi|^2+\E^m[\underset{0\leq t\leq T}{\Sup}|S_{t}|^2]+\Int_{\R^d}\Int_0^T[|f^0(t,x)|^2+|g^0(t,x)|^2+|h^0(t,x)|^2]dtdx\right ]
\end{eqnarray}
\end{Lemma}

The proof of this lemma is a repeat of the proof of Lemma 6 in \cite{MS}. One just needs to notice that when applying the Grownwall's inequality, the constant $c_M$ on
on right of (\ref{2.6}) only depends on the norm of $k$ which is bounded by $M$.
\vskip 0.3cm
We also need the following estimate.
\begin{Lemma}
\begin{equation}\label{2.7}
 \underset{n}{\Sup}\underset{\{k\in K; \|k\|_K\leq M\}}{\Sup}\E^m[n\Int_0^T[(Y_{t}^{n}-S_t)^-]^2dt]\leq C_M.
 \end{equation}
 \end{Lemma}
 {\bf Proof}. Let $F(z)=z^2$. Applying the Ito's formula (see \cite{MS}) we have
 \begin{eqnarray}\label{2.8}
&&F(Y_{t}^{n}-S_t)=F(\Phi(W_{T})-S_T)+\Int_{t}^{T}F^{\prime}(Y_{r}^{n}-S_r)f(r,W_{r},Y_{r}^{n},Z_{r}^{n})dr\nonumber\\
&&+\sum_{j=1}^{\infty}\Int_{t}^{T}F^{\prime}(Y_{r}^{n}-S_r)h_j(r,W_{r},Y_{r}^{n},Z_{r}^{n})k_r^jdr+n\Int_t^T F^{\prime}(Y_{r}^{n}-S_r)(Y^n_r-S_r)^-dr\nonumber\\
&&+\frac{1}{2}\Int_s^T F^{\prime}(Y_{r}^{n}-S_r)g(r,W_{r},Y_{r}^{n},Z_{r}^{n})\ast dW_r\nonumber\\
&&+\Int_t^T<\nabla (F^{\prime}(u^n(r, \cdot)-L(r, \cdot))), g(r,\cdot,u^n(r, \cdot),\nabla u^n(r, \cdot))>(W_r)dr\nonumber\\
&&-\Int_{t}^{T}F^{\prime}(Y_{r}^{n}-S_r)Z_{r}^{n}dW_{r}+\Int_{t}^{T}F^{\prime}(Y_{r}^{n}-S_r)\frac{\partial L}{\partial r}(r,W_r)dr\nonumber\\
&&+\Int_{t}^{T}F^{\prime}(Y_{r}^{n}-S_r)\frac{1}{2}\Delta L(r,W_r)dr+\Int_{t}^{T}F^{\prime}(Y_{r}^{n}-S_r)\nabla L(r,W_r)dW_{r}\nonumber\\
&&-\frac{1}{2}\int_t^TF^{\prime\prime}(Y_{r}^{n}-S_r)|Z_r^n-\nabla L(r,W_r)|^2dr.
\end{eqnarray}
Rearranging the terms we get
 \begin{eqnarray}\label{2.9}
&&(Y_{t}^{n}-S_t)^2 +\int_t^T|Z_r^n-\nabla L(r,W_r)|^2dr+2n\Int_t^T [(Y^n_r-S_r)^-]^2dr\nonumber\\
&&=(\Phi(W_{T})-S_T)^2+2\Int_{t}^{T}(Y_{r}^{n}-S_r)f(r,W_{r},Y_{r}^{n},Z_{r}^{n})dr\nonumber\\
&&+2\sum_{j=1}^{\infty}\Int_{t}^{T}(Y_{r}^{n}-S_r)h_j(r,W_{r},Y_{r}^{n},Z_{r}^{n})k_r^jdr\nonumber\\
&&+\Int_t^T (Y_{r}^{n}-S_r)g(r,W_{r},Y_{r}^{n},Z_{r}^{n})\ast dW_r+2\Int_t^T<Z_r^n-\nabla L(r,W_r), g(r, W_r, Y^n_r, Z_r^n)>dr\nonumber\\
&&-2\Int_{t}^{T}(Y_{r}^{n}-S_r)Z_{r}^{n}dW_{r}+2\Int_{t}^{T}(Y_{r}^{n}-S_r)\frac{\partial L}{\partial r}(r,W_r)dr\nonumber\\
&&+\Int_{t}^{T}(Y_{r}^{n}-S_r)\frac{1}{2}\Delta L(r,W_r)dr+2\Int_{t}^{T}(Y_{r}^{n}-S_r)\nabla L(r,W_r)dW_{r}.\nonumber\\
&&
\end{eqnarray}
Using the conditions on $h$ in the Assumption 2.1, for any given positive constant $\varepsilon_1$ we have
\begin{eqnarray}\label{2.9-1}
&&2\sum_{j=1}^{\infty}\Int_{t}^{T}(Y_{r}^{n}-S_r)h_j(r,W_{r},Y_{r}^{n},Z_{r}^{n})k_r^jdr\nonumber\\
&&=2\Int_{t}^{T}(Y_{r}^{n}-S_r)\sum_{j=1}^{\infty}(h_j(r,W_{r},Y_{r}^{n},Z_{r}^{n})-h_j(r,W_{r},S_r,\nabla L(r,W_r)))k_r^jdr\nonumber\\
&&+2\Int_{t}^{T}(Y_{r}^{n}-S_r)\sum_{j=1}^{\infty}(h_j(r,W_{r},S_r,\nabla L(r,W_r))-h_j(r,W_{r},0,0))k_r^jdr\nonumber\\
&&+2\Int_{t}^{T}(Y_{r}^{n}-S_r)\sum_{j=1}^{\infty}h_j(r,W_{r},0,0)k_r^jdr\nonumber\\
&&\leq 2\Int_{t}^{T}|Y_{r}^{n}-S_r|(\sum_{j=1}^{\infty}(h_j(r,W_{r},Y_{r}^{n},Z_{r}^{n})-h_j(r,W_{r},S_r,\nabla L(r,W_r)))^2)^{\frac{1}{2}}\|k_r\|_{l^2}dr\nonumber\\
&&+2\Int_{t}^{T}|Y_{r}^{n}-S_r|(\sum_{j=1}^{\infty}(h_j(r,W_{r},S_r,\nabla L(r,W_r))-h_j(r,W_{r},0,0))^2)^{\frac{1}{2}}\|k_r\|_{l^2}dr\nonumber\\
&&+2\Int_{t}^{T}|Y_{r}^{n}-S_r|(\sum_{j=1}^{\infty}h_j(r,W_{r},0,0)^2)^{\frac{1}{2}}\|k_r\|_{l^2}dr\nonumber\\
&&\leq C\Int_{t}^{T}|Y_{r}^{n}-S_r|^2\|k_r\|_{l^2}^2dr+\varepsilon_1\Int_{t}^{T}|Z_r^n-\nabla L(r,W_r)|^2dr\nonumber\\
&&+C\Int_{t}^{T}[L(r,W_r)^2+|\nabla L(r,W_r)|^2 +h^0(r,W_{r})^2]dr.
\end{eqnarray}
By the assumptions on $g$, for any given positive constant $\varepsilon_2$ we have
\begin{eqnarray}\label{2.9-2}
&&2\Int_t^T<Z_r^n-\nabla L(r,W_r), g(r, W_r, Y^n_r, Z_r^n)>dr\nonumber\\
&&=2\Int_t^T<Z_r^n-\nabla L(r,W_r), g(r, W_r, Y^n_r, Z_r^n)-g(r,W_{r},S_r,\nabla L(r,W_r))>dr\nonumber\\
&&+2\Int_t^T<Z_r^n-\nabla L(r,W_r), g(r,W_{r},S_r,\nabla L(r,W_r))-g(r,W_{r},0,0>dr\nonumber\\
&&+2\Int_t^T<Z_r^n-\nabla L(r,W_r), g(r,W_{r},0,0>dr\nonumber\\
&&\leq 2C \Int_t^T|Z_r^n-\nabla L(r,W_r)| |Y^n_r-S_r|dr +2\alpha \Int_t^T|Z_r^n-\nabla L(r,W_r)|^2dr\nonumber\\
&&+C\Int_t^T|Z_r^n-\nabla L(r,W_r)|[|L(r,W_r)|+|\nabla L(r,W_r)| +g^0(r,W_{r})]dr\nonumber\\
&&\leq C \Int_t^T|Y^n_r-S_r|^2dr +(2\alpha +\varepsilon_2)\Int_t^T|Z_r^n-\nabla L(r,W_r)|^2dr\nonumber\\
&&+C\Int_t^T[|L(r,W_r)|^2+|\nabla L(r,W_r)|^2 +g^0(r,W_{r})^2]dr.
\end{eqnarray}
By a similar calculation, we have for any given $\varepsilon_3>0$,
\begin{eqnarray}\label{2.9-3}
&&2\Int_{t}^{T}(Y_{r}^{n}-S_r)f(r,W_{r},Y_{r}^{n},Z_{r}^{n})dr\nonumber\\
&&\leq C \Int_t^T|Y^n_r-S_r|^2dr +\varepsilon_3\Int_t^T|Z_r^n-\nabla L(r,W_r)|^2dr\nonumber\\
&&+C\Int_t^T[|L(r,W_r)|^2+|\nabla L(r,W_r)|^2 +f^0(r,W_{r})^2]dr.
\end{eqnarray}
Substitute (\ref{2.9-1}), (\ref{2.9-2}) and (\ref{2.9-3}) back into to (\ref{2.9}), choose $\varepsilon_1, \varepsilon_2, \varepsilon_3$ sufficiently small to obtain
\begin{eqnarray}\label{2.10}
&&\E^m[(Y_{t}^{n}-S_t)^2] +\E^m[\int_t^T|Z_r^n-\nabla L(r,W_r)|^2dr]+2n\E^m[\Int_t^T [(Y^n_r-S_r)^-]^2dr]\nonumber\\
&\leq& C\E^m[(\Phi(W_{T})-S_T)^2]+C\E^m[\Int_{t}^{T}\{f^0(r,W_{r})^2+h^0(r,W_{r})^2+g^0(r,W_{r})^2\}dr\nonumber\\
&+&C\Int_{t}^{T}\E^m[(Y_{r}^{n}-S_r)^2]\|k_r\|_{l^2}^2dr+C\E^m[\Int_{t}^{T}(Y_{r}^{n}-S_r)^2dr]\nonumber\\
&+&C\E^m[\Int_{t}^{T}[(\frac{\partial L}{\partial r}(r,W_r)+\Delta L(r,W_r))^2+|L(r,W_r)|^2+|\nabla L(r,W_r)|^2]dr,
\end{eqnarray}
where the condition on $\alpha$ in the Assumption 2.1 was used.
Now the desired conclusion (\ref{2.7}) follows from the Grownwall's inequality.

\begin{Lemma} For $M>0$, we have
\begin{equation}\label{2.10-1}
 \underset{n \rightarrow \infty}{\lim}\underset{\{k\in K; \|k\|_K\leq M\}}{\Sup}\E^m[\underset{0\leq t\leq T} {\Sup}[(Y_{t}^{n}-S_t)^-]^4] =0.
 \end{equation}
 \end{Lemma}
 {\bf Proof}. Let $G(z)=(z^-)^4$. By the Ito's formula we have
 \begin{eqnarray}\label{2.11}
&&G(Y_{t}^{n}-S_t)=\Int_{t}^{T}G^{\prime}(Y_{r}^{n}-S_r)f(r,W_{r},Y_{r}^{n},Z_{r}^{n})dr\nonumber\\
&&+\sum_{j=1}^{\infty}\Int_{t}^{T}G^{\prime}(Y_{r}^{n}-S_r)h_j(r,W_{r},Y_{r}^{n},Z_{r}^{n})k_r^jdr+n\Int_t^T G^{\prime}(Y_{r}^{n}-S_r)(Y^n_r-S_r)^-dr\nonumber\\
&&+\frac{1}{2}\Int_s^T G^{\prime}(Y_{r}^{n}-S_r)g(r,W_{r},Y_{r}^{n},Z_{r}^{n})\ast dW_r\nonumber\\
&&+\Int_t^T<\nabla (G^{\prime}(u^n(r, \cdot)-L(r, \cdot))), g(r,\cdot,u^{n}(r,\cdot),\nabla u^{n}(r,\cdot))>(W_r)dr\nonumber\\
&&-\Int_{t}^{T}G^{\prime}(Y_{r}^{n}-S_r)Z_{r}^{n}dW_{r}+\Int_{t}^{T}G^{\prime}(Y_{r}^{n}-S_r)\frac{\partial L}{\partial r}(r,W_r)dr\nonumber\\
&&+\Int_{t}^{T}G^{\prime}(Y_{r}^{n}-S_r)\frac{1}{2}\Delta L(r,W_r)dr+\Int_{t}^{T}G^{\prime}(Y_{r}^{n}-S_r)\nabla L(r,W_r)dW_{r}\nonumber\\
&&-\frac{1}{2}\int_t^TG^{\prime\prime}(Y_{r}^{n}-S_r)|Z_r^n-\nabla L(r,W_r)|^2dr.
\end{eqnarray}
Rearrange the terms in the above equation to get
\begin{eqnarray}\label{2.12}
&&[(Y_{t}^{n}-S_t)^-]^4+6\int_t^T[(Y_{r}^{n}-S_r)^-]^2|Z_r^n-\nabla L(r,W_r)|^2dr+4n\Int_t^T [(Y^n_r-S_r)^-]^4dr\nonumber\\
&=&-4\Int_{t}^{T}[(Y_{r}^{n}-S_r)^-]^3f(r,W_{r},Y_{r}^{n},Z_{r}^{n})dr-4\sum_{j=1}^{\infty}\Int_{t}^{T}[(Y_{r}^{n}-S_r)^-]^3h_j(r,W_{r},Y_{r}^{n},Z_{r}^{n})k_r^jdr\nonumber\\
&&-2\Int_t^T [(Y_{r}^{n}-S_r)^-]^3g(r,W_{r}(x),Y_{r}^{n},Z_{r}^{n})\ast dW_r\nonumber\\
&&+12\Int_t^T[(Y_{r}^{n}-S_r)^-]^2<Z^n_r-\nabla L(r,W_r), g(r, W_r, Y^n_r, Z_r^n)>dr\nonumber\\
&&+4\Int_{t}^{T}[(Y_{r}^{n}-S_r)^-]^3Z_{r}^{n}dW_{r}-4\Int_{t}^{T}[(Y_{r}^{n}-S_r)^-]^3\frac{\partial L}{\partial r}(r,W_r)dr\nonumber\\
&&-2\Int_{t}^{T}[(Y_{r}^{n}-S_r)^-]^3\Delta L(r,W_r)dr-4\Int_{t}^{T}[(Y_{r}^{n}-S_r)^-]^3\nabla L(r,W_r)dW_{r}.
\end{eqnarray}
By Assumption 2.1, for any given positive constant $\varepsilon_1$ we have
\begin{eqnarray}\label{2.12-1}
&&12\Int_t^T[(Y_{r}^{n}-S_r)^-]^2<Z_r^n-\nabla L(r,W_r), g(r, W_r, Y^n_r, Z_r^n)>dr\nonumber\\
&&=12\Int_t^T[(Y_{r}^{n}-S_r)^-]^2<Z_r^n-\nabla L(r,W_r), g(r, W_r, Y^n_r, Z_r^n)-g(r,W_{r},S_r,\nabla L(r,W_r))>dr\nonumber\\
&&+12\Int_t^T[(Y_{r}^{n}-S_r)^-]^2<Z_r^n-\nabla L(r,W_r), g(r,W_{r},S_r,\nabla L(r,W_r))-g(r,W_{r},0,0>dr\nonumber\\
&&+12\Int_t^T[(Y_{r}^{n}-S_r)^-]^2<Z_r^n-\nabla L(r,W_r), g(r,W_{r},0,0>dr\nonumber\\
&&\leq C \Int_t^T[(Y_{r}^{n}-S_r)^-]^3|Z_r^n-\nabla L(r,W_r)| dr +12\alpha \Int_t^T[(Y_{r}^{n}-S_r)^-]^2|Z_r^n-\nabla L(r,W_r)|^2dr\nonumber\\
&&+C\Int_t^T[(Y_{r}^{n}-S_r)^-]^2|Z_r^n-\nabla L(r,W_r)|[|L(r,W_r)|+|\nabla L(r,W_r)|]dr\nonumber\\
&&+C \Int_t^T[(Y_{r}^{n}-S_r)^-]^2|Z_r^n-\nabla L(r,W_r)|g^0(r,W_{r})dr\nonumber\\
&&\leq \varepsilon_1\Int_t^T[(Y_{r}^{n}-S_r)^-]^2|Z_r^n-\nabla L(r,W_r)|^2 dr +12\alpha \Int_t^T[(Y_{r}^{n}-S_r)^-]^2|Z_r^n-\nabla L(r,W_r)|^2dr\nonumber\\
&&+C\Int_t^T[(Y_{r}^{n}-S_r)^-]^4[|L(r,W_r)|^2+|\nabla L(r,W_r)|^2+g^0(r,W_{r})^2]dr+C \Int_t^T[(Y_{r}^{n}-S_r)^-]^4dr\nonumber\\
&&\leq (\varepsilon_1+12\alpha)\Int_t^T[(Y_{r}^{n}-S_r)^-]^2|Z_r^n-\nabla L(r,W_r)|^2 dr +C \Int_t^T[(Y_{r}^{n}-S_r)^-]^4dr
\end{eqnarray}
Using again  Assumption 2.1 and the similar computation as above we can show that for any constants $\varepsilon_2>0, \varepsilon_3>0$,
\begin{eqnarray}\label{2.12-2}
&&-4\sum_{j=1}^{\infty}\Int_{t}^{T}[(Y_{r}^{n}-S_r)^-]^3h_j(r,W_{r},Y_{r}^{n},Z_{r}^{n})k_r^jdr\nonumber\\
&&\leq \varepsilon_2\Int_t^T[(Y_{r}^{n}-S_r)^-]^2|Z_r^n-\nabla L(r,W_r)|^2 dr +C \Int_t^T[(Y_{r}^{n}-S_r)^-]^4dr\nonumber\\
&&+C \Int_t^T[(Y_{r}^{n}-S_r)^-]^4\|k_r\|_{l^2}^2dr+C \Int_t^T[(Y_{r}^{n}-S_r)^-]^2dr,
\end{eqnarray}
and
\begin{eqnarray}\label{2.12-3}
&&-4\Int_{t}^{T}[(Y_{r}^{n}-S_r)^-]^3f(r,W_{r},Y_{r}^{n},Z_{r}^{n})dr\nonumber\\
&&\leq \varepsilon_3\Int_t^T[(Y_{r}^{n}-S_r)^-]^2|Z_r^n-\nabla L(r,W_r)|^2 dr +C \Int_t^T[(Y_{r}^{n}-S_r)^-]^4dr\nonumber\\
&&++C \Int_t^T[(Y_{r}^{n}-S_r)^-]^2dr.
\end{eqnarray}

Put (\ref{2.12-3}), (\ref{2.12-2}), (\ref{2.12-1}) and (\ref{2.12}) together, select the constants  $\varepsilon_1$, $\varepsilon_2$ and $\varepsilon_3$ sufficiently small, and take expectation to get
\begin{eqnarray}\label{2.12-4}
&&\E^m[[(Y_{t}^{n}-S_t)^-]^4]+\E^m[\int_t^T[(Y_{r}^{n}-S_r)^-]^2|Z_r^n-\nabla L(r,W_r)|^2dr]+4n\E^m[\Int_t^T [(Y^n_r-S_r)^-]^4dr]\nonumber\\
&&\leq +C \Int_t^T\E^m[[(Y_{r}^{n}-S_r)^-]^4]dr\nonumber\\
&&+C \Int_t^T\E^m[[(Y_{r}^{n}-S_r)^-]^4]\|k_r\|_{l^2}^2dr+C \E^m[\Int_t^T[(Y_{r}^{n}-S_r)^-]^2dr]
\end{eqnarray}
Applying the Grownwall's  inequality and Lemma 4.2 we obtain
\begin{eqnarray}\label{2.13}
 &&\underset{n \rightarrow \infty}{\lim}\underset{\{k\in K; \|k\|_K\leq M\}}{\Sup}\underset{0\leq t\leq T} {\Sup}\E^m[[(Y_{t}^{n}-S_t)^-]^4]\nonumber\\
 &\leq & C_M \underset{n \rightarrow \infty}{\lim}\underset{\{k\in K; \|k\|_K\leq M\}}{\Sup}\E^m[\Int_{0}^{T}[(Y_{r}^{n}-S_r)^-]^2dr]=0,
 \end{eqnarray}
 and
 \begin{equation}\label{2.14}
 \underset{n \rightarrow \infty}{\lim}\underset{\{k\in K; \|k\|_K\leq M\}}{\Sup}\E^m[\int_0^T[(Y_{r}^{n}-S_r)^-]^2|Z_r^n-\nabla L(r,W_r)|^2dr]=0.
\end{equation}
Observe that by the assumptions on the function $g$,
\begin{eqnarray}\label{2.15}
&&2\E^m[ \underset{0\leq t\leq T} {\Sup}|\Int_t^T [(Y_{r}^{n}-S_r)^-]^3g(r,W_{r}(x),Y_{r}^{n},Z_{r}^{n})\ast dW_r|]\nonumber\\
&\leq& C\E^m[\left ( \Int_0^T [(Y_{r}^{n}-S_r)^-]^6|g|^2(r,W_{r}(x),Y_{r}^{n},Z_{r}^{n})dr\right )^{\frac{1}{2}}]\nonumber\\
&\leq& \frac{1}{4}\E^m[ \underset{0\leq r\leq T} {\Sup}[(Y_{r}^{n}-S_r)^-]^4]+C \E^m[\Int_0^T [(Y_{r}^{n}-S_r)^-]^2|g|^2(r,W_{r}(x),Y_{r}^{n},Z_{r}^{n})dr]\nonumber\\
&\leq& \frac{1}{4}\E^m[ \underset{0\leq r\leq T} {\Sup}[(Y_{r}^{n}-S_r)^-]^4]+C \E^m[\Int_0^T [(Y_{r}^{n}-S_r)^-]^4dr]\nonumber\\
&&+C \E^m[\Int_0^T [(Y_{r}^{n}-S_r)^-]^2|Z_{r}^{n}-\nabla L(r,W_r)|^2dr]+C \E^m[\Int_0^T [(Y_{r}^{n}-S_r)^-]^2dr],
\end{eqnarray}
and
\begin{eqnarray}\label{2.16}
&&4\E^m[ \underset{0\leq t\leq T} {\Sup}|\Int_{t}^{T}[(Y_{r}^{n}-S_r)^-]^3<Z_{r}^{n}-\nabla L(r,W_r), dW_{r}>|]\nonumber\\
&\leq& C\E^m[\left ( \Int_0^T [(Y_{r}^{n}-S_r)^-]^6|Z_{r}^{n}-\nabla L(r,W_r)|^2dr\right )^{\frac{1}{2}}]\nonumber\\
&\leq& \frac{1}{4}\E^m[ \underset{0\leq r\leq T} {\Sup}[(Y_{r}^{n}-S_r)^-]^4]+C \E^m[\Int_0^T [(Y_{r}^{n}-S_r)^-]^2|Z_{r}^{n}-\nabla L(r,W_r)|^2dr].
\end{eqnarray}
Using (\ref{2.14})-(\ref{2.16}) and taking supremum over the interval $[0,T]$ in (\ref{2.12}) we further deduce that
$$
 \underset{n \rightarrow \infty}{\lim}\underset{\{k\in K; \|k\|_K\leq M\}}{\Sup}\E^m[\underset{0\leq t\leq T} {\Sup}[(Y_{t}^{n}-S_t)^-]^4] =0.
$$
completing the proof.
\begin{Proposition}
 For any $M>0$, we have
 \begin{equation}\label{2.17}
 \underset{n \rightarrow \infty}{\lim}\underset{\{k\in K; \|k\|_K\leq M\}}{\Sup}|u^n-u|_{H_T}=0.
 \end{equation}
 \end{Proposition}
 \noindent{\bf Proof}. We note that
 \begin{eqnarray}\label{2.17-1}
 &&|u^n-u^q|_{H_T}\nonumber\\
 &\leq&\E^m[ \underset{0\leq r\leq T} {\Sup}(Y_{r}^{n}-Y_r^q)^2]+C \E^m[\Int_0^T |Z_r^n-Z_{r}^{q}|^2dr].
 \end{eqnarray}
 In order to prove (\ref{2.17}), it is sufficient to show
  \begin{equation}\label{2.17-2}
  \underset{n, q \rightarrow \infty}{\lim}\underset{\{k\in K; \|k\|_K\leq M\}}{\Sup}\E^m[\underset{0\leq t\leq T} {\Sup}(Y_{t}^{n}-Y_t^q)^2]] =0,
\end{equation}
and
\begin{equation}\label{2.17-3}
  \underset{n, q \rightarrow \infty}{\lim}\underset{\{k\in K; \|k\|_K\leq M\}}{\Sup}\E^m[\Int_0^T |Z_r^n-Z_{r}^{q}|^2dr]=0.
\end{equation}

 We will achieve this with the help of backward stochastic differential equations satisfied by $Y^n_t=u^n(t, W_t)$. Applying Ito's formula we have
\begin{eqnarray}\label{2.18}
&&(Y_{t}^{n}-Y_{t}^{q})^2+\int_t^T|Z_r^n-Z_r^q|^2dr\nonumber\\
&=&2\Int_{t}^{T}(Y_{r}^{n}-Y_r^q)(f(r,W_{r},Y_{r}^{n},Z_{r}^{n})-f(r,W_{r},Y_{r}^{q},Z_{r}^{q}))dr\nonumber\\
&&+2\sum_{j=1}^{\infty}\Int_{t}^{T}(Y_{r}^{n}-Y_r^q)(h_j(r,W_{r},Y_{r}^{n},Z_{r}^{n})-h_j(r,W_{r},Y_{r}^{q},Z_{r}^{q}))k_r^jdr\nonumber\\
&&+2n\Int_t^T (Y_{r}^{n}-Y_r^q)(Y^n_r-S_r)^-dr-2q\Int_t^T (Y_{r}^{n}-Y_r^q)(Y^q_r-S_r)^-dr\nonumber\\
&&+\Int_t^T (Y_{r}^{n}-Y_r^q)(g(r,W_{r},Y_{r}^{n},Z_{r}^{n})-g(r,W_{r},Y_{r}^{m},Z_{r}^{q}))\ast dW_r\nonumber\\
&&+2\Int_t^T<Z_r^n-Z_r^q, g(r,W_{r},Y_{r}^{n},Z_{r}^{n})-g(r,W_{r},Y_{r}^{q},Z_{r}^{q})>dr\nonumber\\
&&-2\Int_{t}^{T}(Y_{r}^{n}-Y_r^q)<Z_{r}^{n}-Z_r^q,dW_{r}>\nonumber\\
&:=& I_1^{n,q}(t)+I_2^{n,q}(t)+I_3^{n,q}(t)+I_4^{n,q}(t)+I_5^{n,q}(t)+I_6^{n,q}(t)+I_7^{n,q}(t).
\end{eqnarray}
Note that
\begin{eqnarray}\label{2.19}
&&I_3^{n,q}(t)+I_4^{n,q}(t)\nonumber\\
&=&2n\Int_t^T (Y_{r}^{n}-Y_r^q)(Y^n_r-S_r)^-dr-2q\Int_t^T (Y_{r}^{n}-Y_r^q)(Y^q_r-S_r)^-dr\nonumber\\
&\leq &2n\Int_t^T (Y_r^q-S_r)^-(Y^n_r-S_r)^-dr+2q\Int_t^T (Y_{r}^{n}-S_r)^-(Y^q_r-S_r)^-dr\nonumber\\
&\leq &2\underset{0\leq r\leq T}{\Sup}(Y_r^q-S_r)^- n\Int_0^T(Y^n_r-S_r)^-dr+2\underset{0\leq r\leq T}{\Sup}(Y_r^n-S_r)^- q\Int_0^T(Y^q_r-S_r)^-dr.
\end{eqnarray}
By Young's inequality, we have for any $\delta_1>0$,
\begin{equation}\label{2.20}
I_1^{n,q}(t)\leq \delta_1\Int_{t}^{T}|Z_{r}^{n}-Z_r^q|^2dr+ C\Int_{t}^{T}|Y_{r}^{n}-Y_r^q|^2dr.
\end{equation}
Moreover for any $\delta_2>0$, we have
\begin{equation}\label{2.21}
I_2^{n,q}(t)\leq \delta_2\Int_{t}^{T}|Z_{r}^{n}-Z_r^q|^2dr+ C\Int_{t}^{T}|Y_{r}^{n}-Y_r^q|^2(1+\|k_r\|_{l^2}^2)dr.
\end{equation}
Using Young's inequality again, we have for any $\delta_3>0$,
\begin{equation}\label{2.22}
I_6^{n,q}(t)\leq (\delta_3+2\alpha)\Int_{t}^{T}|Z_{r}^{n}-Z_r^q|^2dr+ C\Int_{t}^{T}|Y_{r}^{n}-Y_r^q|^2dr.
\end{equation}
Substitute (\ref{2.19})-(\ref{2.22}) back to (\ref{2.18}), choose constants $\delta_i, i=1,2,3 $ sufficiently small and take expectation to obtain

\begin{eqnarray}\label{2.22-1}
&&\E^m[(Y_{t}^{n}-Y_{t}^{q})^2]+\E^m[\int_t^T|Z_r^n-Z_r^q|^2dr]\nonumber\\
&\leq &C\E^m[\Int_{t}^{T}|Y_{r}^{n}-Y_r^q|^2(1+\|k_r\|_{l^2}^2)dr]\nonumber\\
&+& C(\E^m[\underset{0\leq r\leq T}{\Sup}[(Y_r^q-S_r)^-]^2])^{\frac{1}{2}} (\E^m[(n\Int_0^T(Y^n_r-S_r)^-dr)^2])^{\frac{1}{2}}\nonumber\\
&+&C(\E^m[\underset{0\leq r\leq T}{\Sup}[(Y_r^n-S_r)^-]^2])^{\frac{1}{2}} (\E^m[(q\Int_0^T(Y^q_r-S_r)^-dr)^2])^{\frac{1}{2}}
\end{eqnarray}
Using Lemma 4.1, Lemma 4.3  and  applying the Grownwall's inequality we deduce that
\begin{equation}\label{2.22-2}
 \underset{n, q \rightarrow \infty}{\lim}\underset{\{k\in K; \|k\|_K\leq M\}}{\Sup}\underset{0\leq t\leq T} {\Sup}\E^m[(Y_{t}^{n}-Y_t^q)^2]] =0,
\end{equation}
and
\begin{equation}\label{2.22-3}
 \underset{n, q \rightarrow \infty}{\lim}\underset{\{k\in K; \|k\|_K\leq M\}}{\Sup}\E^m[\int_0^T|Z_{t}^{n}-Z_t^q|^2dt] =0.
\end{equation}
Next we will strengthen the convergence in (\ref{2.22-2}) to
\begin{equation}\label{2.22-4}
 \underset{n, q \rightarrow \infty}{\lim}\underset{\{k\in K; \|k\|_K\leq M\}}{\Sup}\E^m[\underset{0\leq t\leq T} {\Sup}(Y_{t}^{n}-Y^q)^2]=0.
\end{equation}
We notice that by the Burkh\"o{}lder's inequality, for any $\delta_4>0$ we have
\begin{eqnarray}\label{2.23}
&&\E^m[ \underset{0\leq t\leq T} {\Sup}|I_5^{n,q}(t)|]\nonumber\\
&\leq &C\E^m[\left ( \Int_0^T (Y_{r}^{n}-Y_r^q)^2|g(r,W_r,Y_{r}^{n},Z_{r}^{n})-g(r,W_{r},Y_{r}^{q},Z_{r}^{q})|^2dr\right )^{\frac{1}{2}}]\nonumber\\
&\leq& \delta_4\E^m[ \underset{0\leq r\leq T} {\Sup}(Y_{r}^{n}-Y_r^q)^2]+C \E^m[\Int_0^T |g(r,W_{r},Y_{r}^{n},Z_{r}^{n})-g(r,W_{r},Y_{r}^{q},Z_{r}^{q})|^2dr]\nonumber\\
&\leq&\delta_4\E^m[ \underset{0\leq r\leq T} {\Sup}(Y_{r}^{n}-Y_r^q)^2]+C\E^m[\Int_0^T |Z_{r}^{n}-Z_{r}^{q})|^2dr]\nonumber\\
&&+C\E^m[\Int_0^T |Y_{r}^{n}-Y_{r}^{q})|^2dr].
\end{eqnarray}
Similarly, we have for $\delta_5>0$
\begin{eqnarray}\label{2.24}
&&\E^m[ \underset{0\leq t\leq T} {\Sup}|I_7^{n,q}(t)|]\nonumber\\
&\leq &C\E^m[\left ( \Int_0^T (Y_{r}^{n}-Y_r^q)^2|Z_r^n-Z_r^q|^2dr\right )^{\frac{1}{2}}]\nonumber\\
&\leq& \delta_5\E^m[ \underset{0\leq r\leq T} {\Sup}(Y_{r}^{n}-Y_r^q)^2]+C \E^m[\Int_0^T |Z_r^n-Z_{r}^{q}|^2dr].
\end{eqnarray}
Now use the above two estimates (\ref{2.23}) and (\ref{2.24}) and the already proved (\ref{2.22-2}) to obtain (\ref{2.22-4}).
This completes the proof.
\begin{Theorem}
Let  Assumptions 2.1 hold. Assume that $k^{\varepsilon}\rightarrow k$ weakly in the Hilbert space $K$ as $\varepsilon \rightarrow 0$.
Then $u^{\varepsilon}$ converges to $u$ in the space $H_T$, where $u^{\varepsilon}$ denotes the solution of equation (\ref{2.1}) with $k^{\varepsilon}$ replacing $k$.
\end{Theorem}
 \noindent{\bf Proof}. We will first prove a similar convergence result for the corresponding penalized PDEs and then combined with the uniform convergence proved in Proposition 4.1 we complete the proof of Theorem 4.1.   Let $u^{\varepsilon, n}$ be the solution to the following penalized PDE:
 \begin{eqnarray}\label{2.25}
du^{\varepsilon, n} (t,x)&+&\frac{1}{2}\Delta u^{\varepsilon, n}(t,x)+\sum_{i=1}^d\partial_ig_i(t,x, u^{\varepsilon, n}(t,x), \nabla u^{\varepsilon, n}(t,x))dt +f(t,x,u^{\varepsilon, n}(t,x), \nabla u^{\varepsilon, n}(t,x))dt\nonumber\\
 &+ &\sum_{j=1}^{\infty} h_j(t,x,u^{\varepsilon, n}(t,x), \nabla u^{\varepsilon, n}(t,x))k^{\varepsilon,j}_tdt =-n(u^{\varepsilon, n}(t,x)-L(t,x))^-dt ,\\
 u^{\varepsilon, n}(T,x)&=&\Phi (x), \quad\quad x\in \R^d.
 \end{eqnarray}
 We first fix the integer $n$ and show $\lim_{\varepsilon\rightarrow \infty}\|u^{\varepsilon,n}-u^n\|_{H_T}=0$, $u^n$ is the solution of equation (\ref{2.25}) with $k^{\varepsilon}$ replaced by $k$.  To this end, we first prove that the family $\{ u^{\varepsilon, n}, \varepsilon>0\}$ is tight in the space $L^2([0, T], L_{loc}^2(\R^d))$. Using the chain rule and Gronwall's inequality, as in Lemma 4.1 , we can show that
 \begin{equation}\label{2.26}
 \underset{\varepsilon} {\Sup}\|u^{\varepsilon, n}\|_{H_T}^2=\underset{\varepsilon} {\Sup}\{\sup_{0\leq t\leq T}|u^{\varepsilon, n}(t)|^2+\int_0^T\|u^{\varepsilon, n}(t)\|^2dt\}<\infty.
 \end{equation}
 For $\beta\in (0,1)$, recall that $W^{\beta,2}([0,T], V^*)$ is the space of mappings $v(\cdot): [0, T]\rightarrow V^*$ that satisfy
 \begin{equation}\label{2.27}
 \|v\|_{W^{\beta,2}([0,T], V^*)}^2=\int_0^T\|v(t)\|_{V^*}^2+\int_0^T\int_0^T\frac{\|v(t)-v(s)\|_{V^*}^2}{|t-s|^{1+2\beta}}  <\infty.
 \end{equation}
 It is well known  (see e.g. \cite{FG}) that the imbedding
 $$ L^2([0,T], V)\cap W^{\beta,2}([0,T], V^*)\hookrightarrow L^2([0,T], L_{loc}^2(\R^d))$$
 is compact.
 As an equation in $V^*$, we have
 \begin{eqnarray}\label{2.28}
u^{\varepsilon, n} (t)&=&\Phi+\frac{1}{2}\int_t^T\Delta u^{\varepsilon, n}(s)ds+\int_t^T\sum_{i=1}^d\partial_ig_i(s,x, u^{\varepsilon, n}(s,x), \nabla u^{\varepsilon, n}(s,x))ds \nonumber\\
&+&\int_t^Tf(s,x,u^{\varepsilon, n}(s,x), \nabla u^{\varepsilon, n}(s,x))ds\nonumber\\
 &+ &\sum_{j=1}^{\infty} \int_t^Th_j(s,x,u^{\varepsilon, n}(s,x), \nabla u^{\varepsilon, n}(s,x))k^{\varepsilon,j}_sds +n\int_t^T(u^{\varepsilon, n}(s,x)-L(s,x))^-ds \nonumber\\
 &:=& \Phi+I_1(t)+I_2(t)+I_3(t)+I_4(t)+I_5(t).
 \end{eqnarray}
 In view of (\ref{2.26}), we have
 \begin{equation}\label{2.29}
 \|I_1(t)-I_1(s)\|_{V^*}^2\leq C\int_s^t\|\Delta u^{\varepsilon, n}(r)\|_{V^*}^2dr (|t-s|)\leq C\int_0^T\|u^{\varepsilon, n}(r)\|^2dr (|t-s|)\leq C|t-s|.
 \end{equation}
 Using the condition (iii) in Assumption 2.1, we have
 \begin{equation}\label{2.30}
 \|I_4(t)-I_4(s)\|_{V^*}^2\leq C(\int_0^T\|k^{\varepsilon}_r\|_{l^2}^2dr)|\bar{h}|^2 |t-s|\leq C|t-s|.
 \end{equation}
 By (\ref{2.26}) and the similar calculations as above  we also have
 \begin{equation}\label{2.31}
 \|I_i(t)-I_i(s)\|_{V^*}^2\leq C|t-s|, \quad\quad\quad\quad \quad i=2,3,5.
 \end{equation}
 Thus, for $\beta\in (0,\frac{1}{2})$, it follows from (\ref{2.28}) --(\ref{2.31})  that
 \begin{equation}\label{2.32}
 \underset{\varepsilon} {\Sup}\|u^{\varepsilon, n}\|_{W^{\beta,2}([0,T], V^*)}^2<\infty.
 \end{equation}
 Combining (\ref{2.32}) with (\ref{2.26}), we conclude that  $\{ u^{\varepsilon, n}, \varepsilon>0\}$ is tight in the space $L^2([0, T], L_{loc}^2(\R^d))$. Now, applying the chain rule, we obtain
 \begin{eqnarray}\label{2.33}
&&|u^{\varepsilon, n}(t)-u^{n}(t)|^2\nonumber\\
&=&-\int_t^T|\nabla (u^{\varepsilon, n}(s)-u^{n}(s))|^2ds\nonumber\\
&&-2\int_t^T<g(s,\cdot, u^{\varepsilon, n}(s,\cdot), \nabla u^{\varepsilon, n}(s,\cdot))-g(s,\cdot, u^{n}(s,\cdot), \nabla u^{n}(s,\cdot)), \nabla (u^{\varepsilon, n}(s)-u^{n}(s))>ds\nonumber\\
 &&+2\int_t^T<f(s,\cdot,u^{\varepsilon, n}(s,\cdot), \nabla u^{\varepsilon, n}(s,\cdot))-f(s,\cdot,u^{n}(s,\cdot), \nabla u^{n}(s,\cdot)),u^{\varepsilon, n}(s)-u^{n}(s)> ds\nonumber\\
 &&+2\int_t^T<u^{\varepsilon, n}(s)-u^{n}(s), \sum_{j=1}^{\infty}(h_j(s,\cdot,u^{\varepsilon, n}(s,\cdot), \nabla u^{\varepsilon, n}(s,\cdot))-h_j(s,\cdot,u^{n}(s,\cdot), \nabla u^{n}(s,\cdot)))k^{\varepsilon,j}_s>ds \nonumber\\
 &&+2\int_t^T<u^{\varepsilon, n}(s)-u^{n}(s), \sum_{j=1}^{\infty}h_j(s,\cdot,u^{n}(s,\cdot), \nabla u^{n}(s,\cdot))(k^{\varepsilon,j}_s-k^{j}_s)>ds  \nonumber\\  &&+2n\int_t^T<u^{\varepsilon, n}(s)-u^{n}(s),   (u^{\varepsilon, n}(s)-L(s,\cdot))^--(u^{n}(s,\cdot)-L(s, \cdot))^->ds
 \end{eqnarray}
 By the assumptions on $h_j$ and Young's inequality, we see that for any given $\delta_1>0$,
 \begin{eqnarray}\label{2.34}
 &&2\int_t^T<u^{\varepsilon, n}(s)-u^{n}(s), \sum_{j=1}^{\infty}(h_j(s,\cdot,u^{\varepsilon, n}(s,\cdot), \nabla u^{\varepsilon, n}(s,\cdot))-h_j(s,\cdot,u^{n}(s,\cdot), \nabla u^{n}(s,\cdot)))k^{\varepsilon,j}_s>ds \nonumber\\
 &\leq& \delta_1 \int_t^T |\nabla (u^{\varepsilon, n}(s)-u^{n}(s))|^2ds+C\int_t^T|u^{\varepsilon, n}(s)-u^{n}(s)|^2(1+\|k_s^{\varepsilon}\|_{l^2}^2)ds.
 \end{eqnarray}
 Using the assumptions on $f$, $g$ and (\ref{2.34}) it follows  from (\ref{2.33}) that there exist positive constants $\delta$, $C$ such that
 \begin{eqnarray}\label{2.35}
&&|u^{\varepsilon, n}(t)-u^{n}(t)|^2+\delta \int_t^T|\nabla (u^{\varepsilon, n}(s)-u^{n}(s))|^2ds\nonumber\\
&\leq& C\int_t^T|u^{\varepsilon, n}(s)-u^{n}(s)|^2(1+\|k_s^{\varepsilon}\|_{l^2}^2)ds\nonumber\\
 &&+2\int_t^T<u^{\varepsilon, n}(s)-u^{n}(s), \sum_{j=1}^{\infty}h_j(s,\cdot,u^{n}(s,\cdot), \nabla u^{n}(s,\cdot))(k^{\varepsilon,j}_s-k^{j}_s)>ds.
 \end{eqnarray}
 By Gronwall's inequality, (\ref{2.35}) yields that
 \begin{eqnarray}\label{2.36}
&&\underset{0\leq t\leq T} {\Sup}|u^{\varepsilon, n}(t)-u^{n}(t)|^2\nonumber\\
&\leq& exp(C\int_0^T(1+\|k_s^{\varepsilon}\|_{l^2}^2)ds)\underset{0\leq t\leq T} {\Sup}|\int_t^T<u^{\varepsilon, n}(s)-u^{n}(s), \sum_{j=1}^{\infty}h_j(s,\cdot,u^{n}(s,\cdot), \nabla u^{n}(s,\cdot))(k^{\varepsilon,j}_s-k^{j}_s)>ds|\nonumber\\
&\leq& C\underset{0\leq t\leq T} {\Sup}|\int_t^T<u^{\varepsilon, n}(s)-u^{n}(s), \sum_{j=1}^{\infty}h_j(s,\cdot,u^{n}(s,\cdot), \nabla u^{n}(s,\cdot))(k^{\varepsilon,j}_s-k^{j}_s)>ds|.
 \end{eqnarray}
 To show $\lim_{\varepsilon\rightarrow 0}\|u^{\varepsilon,n}-u^n\|_{H_T}=0$, in view of (\ref{2.35}) and (\ref{2.36}), it suffices to prove
 \begin{equation}\label{2.37}
 \lim_{\varepsilon\rightarrow 0}\underset{0\leq t\leq T} {\Sup}|\int_t^T<u^{\varepsilon, n}(s)-u^{n}(s), \sum_{j=1}^{\infty}h_j(s,\cdot,u^{n}(s,\cdot), \nabla u^{n}(s,\cdot))(k^{\varepsilon,j}_s-k^{j}_s)>ds|=0.
 \end{equation}
 This will be achieved if we show that for any sequence $\varepsilon_m\rightarrow 0$, one can find a subsequence $\varepsilon_{m_k}\rightarrow 0$ such that
\begin{equation}\label{2.38}
 \lim_{k\rightarrow \infty}\underset{0\leq t\leq T} {\Sup}|\int_t^T<u^{\varepsilon_{m_k} , n}(s)-u^{n}(s), \sum_{j=1}^{\infty}h_j(s,\cdot,u^{n}(s,\cdot), \nabla u^{n}(s,\cdot))(k^{\varepsilon_{m_k},j}_s-k^{j}_s)>ds|=0
 \end{equation}
 Now fix a sequence $\varepsilon_m\rightarrow 0$. Since $\{ u^{\varepsilon_{m}, n}, m\geq 1\}$ is tight in $L^2([0,T], L_{loc}^2(\R^d))$, there exist a subsequence $m_k, k\geq 1$ and a mapping $\tilde{u}$ such that $u^{\varepsilon_{m_k},n}\rightarrow \tilde{u}$ in $L^2([0,T], L_{loc}^2(\R^d))$. Moreover, because of the uniform bound of $u^{\varepsilon_{m_k},n}$ in (\ref{2.26}), $\tilde{u}$ belongs to  $L^2([0,T], H)$.  Now,
 \begin{eqnarray}\label{2.39}
&&\underset{0\leq t\leq T} {\Sup}|\int_t^T<u^{\varepsilon_{m_k} , n}(s)-u^{n}(s), \sum_{j=1}^{\infty}h_j(s,\cdot,u^{n}(s,\cdot), \nabla u^{n}(s,\cdot))(k^{\varepsilon_{m_k},j}_s-k^{j}_s)>ds\nonumber\\
&\leq& \underset{0\leq t\leq T} {\Sup}|\int_t^T<u^{\varepsilon_{m_k} , n}(s)-\tilde{u}(s), \sum_{j=1}^{\infty}h_j(s,\cdot,u^{n}(s,\cdot), \nabla u^{n}(s,\cdot))(k^{\varepsilon_{m_k},j}_s-k^{j}_s)>ds\nonumber\\
&+&\underset{0\leq t\leq T} {\Sup}|\int_t^T<\tilde{u}(s)-u^{n}(s), \sum_{j=1}^{\infty}h_j(s,\cdot,u^{n}(s,\cdot), \nabla u^{n}(s,\cdot))(k^{\varepsilon_{m_k},j}_s-k^{j}_s)>ds.
 \end{eqnarray}
 Since $k^{\varepsilon_{m_k}}\rightarrow k$ weakly in $L^2([0,T], l^2)$, for every $t>0$, it holds that
\begin{equation}\label{2.41}
\lim_{k\rightarrow \infty}\int_t^T<\tilde{u}(s)-u^{n}(s), \sum_{j=1}^{\infty}h_j(s,\cdot,u^{n}(s,\cdot), \nabla u^{n}(s,\cdot))(k^{\varepsilon_{m_k},j}_s-k^{j}_s)>ds=0.
\end{equation}
On the other hand, using the assumption on $h$, for $0<t_1<t_2\leq T$, we have
\begin{eqnarray}\label{2.42}
&&|\int_{t_1}^{t_2}<\tilde{u}(s)-u^{n}(s), \sum_{j=1}^{\infty}h_j(s,\cdot,u^{n}(s,\cdot), \nabla u^{n}(s,\cdot))(k^{\varepsilon_{m_k},j}_s-k^{j}_s)>ds|\nonumber\\
&\leq& C(\int_{t_1}^{t_2}|\tilde{u}(s)-u^{n}(s)|^2ds)^{\frac{1}{2}}(\int_{t_1}^{t_2}\|k^{\varepsilon_{m_k}}_s-k_s\|_{l^2}^2ds)^{\frac{1}{2}}\leq C(\int_{t_1}^{t_2}|\tilde{u}(s)-u^{n}(s)|^2ds)^{\frac{1}{2}}.
\end{eqnarray}
Combing (\ref{2.41}) and (\ref{2.42}) we deduce that
\begin{equation}\label{2.43}
\lim_{k\rightarrow \infty}\underset{0\leq t\leq T} {\Sup}|\int_t^T<\tilde{u}(s)-u^{n}(s), \sum_{j=1}^{\infty}h_j(s,\cdot,u^{n}(s,\cdot), \nabla u^{n}(s,\cdot))(k^{\varepsilon_{m_k},j}_s-k^{j}_s)>ds|=0.
\end{equation}
By H\"o{}lder's inequality and the assumption on $h$, we have
\begin{eqnarray}\label{2.43-1}
&& |\int_t^T<u^{\varepsilon_{m_k} , n}(s)-\tilde{u}(s), \sum_{j=1}^{\infty}h_j(s,\cdot,u^{n}(s,\cdot), \nabla u^{n}(s,\cdot))(k^{\varepsilon_{m_k},j}_s-k^{j}_s)>ds|\nonumber\\
&\leq&\int_0^T\int_{\R^d}|u^{\varepsilon_{m_k} , n}(s,x)-\tilde{u}(s,x)| (\sum_{j=1}^{\infty}h_j^2(s,\cdot,u^{n}(s,x), \nabla u^{n}(s,x))^2)^{\frac{1}{2}}dx (\|k^{\varepsilon_{m_k}}_s\|_{l^2}+\|k_s\|_{l^2})ds)\nonumber\\
&\leq& (\int_0^T(\|k^{\varepsilon_{m_k}}_s\|_{l^2}^2+\|k_s\|_{l^2}^2)ds)^{\frac{1}{2}}
(\int_0^Tds(\int_{\R^d}|u^{\varepsilon_{m_k} , n}(s,x)-\tilde{u}(s,x)| \bar{h}(x)dx)^2)^{\frac{1}{2}}\nonumber\\
&\leq&C (\int_0^Tds(\int_{\R^d}|u^{\varepsilon_{m_k} , n}(s,x)-\tilde{u}(s,x)| \bar{h}(x)dx)^2)^{\frac{1}{2}}.
 \end{eqnarray}
For any $M>0$,  denote by $B_M$ the ball in $\R^d$ centered at zero with radius $M$.  We can bound the right side of (\ref{2.43-1}) as follows:
 \begin{eqnarray}\label{2.43-2}
&&\int_0^Tds(\int_{\R^d}|u^{\varepsilon_{m_k} , n}(s,x)-\tilde{u}(s,x)| \bar{h}(x)dx)^2\nonumber\\
&\leq& C\int_0^Tds(\int_{B_M}|u^{\varepsilon_{m_k} , n}(s,x)-\tilde{u}(s,x)|^2dx) (\int_{\R^d}\bar{h}^2(x)dx)\nonumber\\
&&+C\int_0^Tds(\int_{\R^d}|u^{\varepsilon_{m_k} , n}(s,x)-\tilde{u}(s,x)|^2dx) (\int_{B_M^c}\bar{h}^2(x)dx\nonumber\\
&\leq& C\int_0^Tds(\int_{B_M}|u^{\varepsilon_{m_k} , n}(s,x)-\tilde{u}(s,x)|^2dx)+C\int_{B_M^c}\bar{h}^2(x)dx,
 \end{eqnarray}
 where the uniform $L^2([0, T]\times \R^d)$-bound of $u^{\varepsilon_{m_k} , n}$ has been used. Now given any constant $\delta>0$, we can pick a constant $M$ such that $C\int_{B_M^c}\bar{h}^2(x)dx\leq \delta$. For the chosen constant $M$, we have
 $$\lim_{k\rightarrow \infty}\int_0^Tds(\int_{B_M}|u^{\varepsilon_{m_k} , n}(s,x)-\tilde{u}(s,x)|^2dx=0.$$
 Thus, it follows from (\ref{2.43-1}), (\ref{2.43-2}) that
 \begin{equation}\label{2.43-3}
\lim_{k\rightarrow \infty}\underset{0\leq t\leq T} {\Sup}|\int_t^T<u^{\varepsilon_{m_k} , n}(s)-\tilde{u}(s), \sum_{j=1}^{\infty}h_j(s,\cdot,u^{n}(s,\cdot), \nabla u^{n}(s,\cdot))(k^{\varepsilon_{m_k},j}_s-k^{j}_s)>ds\leq \delta^{\frac{1}{2}}.
\end{equation}
Since $\delta$ is arbitrary, (\ref{2.38}) follows from (\ref{2.39}), (\ref{2.43}) and (\ref{2.43-3}). Hence we have proved $\lim_{\varepsilon\rightarrow 0}\|u^{\varepsilon,n}-u^n\|_{H_T}=0$.
 \vskip 0.2cm
Now we are ready to complete the last step of the proof.  For any $n\geq 1$, we have
\begin{eqnarray}\label{2.44}
&&\|u^{\varepsilon}-u\|_{H_T}\nonumber\\
&\leq&\|u^{\varepsilon}-u^{\varepsilon, n}\|_{H_T}+\|u^{\varepsilon, n}-u^{n}\|_{H_T}+\|u^{n}-u\|_{H_T}.
\end{eqnarray}
For any given $\delta>0$, by Proposition 4.1 there exists an integer $n_0$ such that $\sup_{\varepsilon}\|u^{\varepsilon}-u^{\varepsilon, n_0}\|_{H_T}\leq \frac{\delta}{2}$ and $\|u-u^{n_0}\|_{H_T}\leq \frac{\delta}{2}$. Replacing $n$ in (\ref{2.44}) by $n_0$ we get
$$\|u^{\varepsilon}-u\|_{H_T}\leq \delta+\|u^{\varepsilon, n_0}-u^{n_0}\|_{H_T}.$$ As we just proved
$$\lim_{\varepsilon \rightarrow 0}\|u^{\varepsilon, n_0}-u^{n_0}\|_{H_T}=0,$$
we obtain that
$$\lim_{\varepsilon \rightarrow 0}\|u^{\varepsilon}-u\|_{H_T}\leq \delta.$$
Since the constant $\delta$ is arbitrary, the proof is complete.

\section{Large deviations}
After the preparations in Section 4, we are ready to state  and to prove the large deviation result. Recall that $U^{\varepsilon}$ is the solution of the obstacle problem:
\begin{eqnarray}\label{3.1-1}
dU^{\varepsilon}(t,x)&+&\frac{1}{2}\Delta U^{\varepsilon}(t,x)+\sum_{i=1}^d\partial_ig_i(t,x,U^{\varepsilon}(t,x),\nabla U^{\varepsilon}(t,x))dt +f(t,x, U^{\varepsilon}(t,x),\nabla U^{\varepsilon}(t,x))dt\nonumber\\
 &+ &\sqrt{\varepsilon}\sum_{j=1}^{\infty} h_j(t,x, U^{\varepsilon}(t,x),\nabla U^{\varepsilon}(t,x))dB^j_t =-R^{\varepsilon}(dt,dx),\\
 U^{\varepsilon}(t,x)&\geq & L(t,x), \quad \quad (t,x)\in \R^+\times \R^d,\nonumber\\
 U^{\varepsilon}(T,x)&=&\Phi (x), \quad\quad x\in \R^d.
 \end{eqnarray}
For $k\in K=L^2([0,T], l^2)$, denote by $u^k$ the solution of the following  deterministic obstacle problem:
\begin{eqnarray}\label{3.1-2}
du^k(t,x)&+&\frac{1}{2}\Delta u^k(t,x)+\sum_{i=1}^d\partial_ig_i(t,x, u^k(t,x), \nabla u^k(t,x))dt +f(t,x,u^k(t,x),\nabla u^k(t,x))dt\nonumber\\
 &+ &\sum_{j=1}^{\infty} h_j(t,x,u^k(t,x),\nabla u^k(t,x))k^j_tdt =-\nu^k(dt,dx),\\
 u^k(t,x)&\geq & L(t,x), \quad \quad (t,x)\in \R^+\times \R^d,\nonumber\\
 u^k(T,x)&=&\Phi (x), \quad\quad x\in \R^d.
 \end{eqnarray}
Define a measurable mapping $\Gamma^0:C([0,T];\mathbb{R}^\infty)\rightarrow H_T$ by
$$
\Gamma^0(\int_0^{\cdot}k_sds):=u^k\quad\quad\quad\quad \mbox{for}\quad\quad\quad k\in K,$$
where $u^k$ is the solution of (\ref{3.1-2}). Here is the main result:
\begin{Theorem}
Let the Assumption 2.1 hold.  Then the family $\{U^\varepsilon\}_{\varepsilon>0}$ satisfies a large deviation principle on the space $H_T$ with the rate function $I$ given by
\begin{equation}
I(g):=\inf_{\{k\in K;g=\Gamma^0(\int_0^{\cdot}k_sds)\}}\left\{\frac12\int_0^T\|k_s\|_{l^2}^2ds\right\},\ g\in H_T,
\end{equation}
with the convention $\inf\{\emptyset\}=\infty$.
 \end{Theorem}
\noindent{\bf Proof}. The existence of a unique strong solution of the obstacle problem (\ref{3.1-1}) means that for every $\varepsilon>0$, there exists
a measurable mapping $\Gamma^\varepsilon\left(\cdot\right):C([0,T];\mathbb{R}^\infty)\rightarrow H_T$ such that
$$U^{\varepsilon}=\Gamma^\varepsilon\left(B(\cdot)\right).$$
To prove the theorem, we are going to show that the conditions (i) and (ii) in Theorem 3.2 are satisfied. Condition (ii) is exactly the statement of Theorem 4.1.
It remains to establish the condition (i) in Theorem 3.2. Recall the definitions of the spaces $S_N$ and $\tilde{S}_N$ given in Section 3. Let $\{k^{\varepsilon}, \varepsilon>0\}\subset \tilde{S}_N$ be a given family of stochastic processes. Applying  Girsanov theorem it is easy to see that $u^{\varepsilon}=\Gamma^\varepsilon\left(B(\cdot)+\frac{1}{\sqrt\varepsilon}\int_0^{\cdot}k^\varepsilon(s)ds\right)$
is the solution of the stochastic obstacle problem:
\begin{eqnarray}\label{3.1}
du^{\varepsilon}(t,x)&+&\frac{1}{2}\Delta u^{\varepsilon}(t,x)+\sum_{i=1}^d\partial_ig_i(t,x,u^{\varepsilon}(t,x),\nabla u^{\varepsilon}(t,x))dt +f(t,x, u^{\varepsilon}(t,x),\nabla u^{\varepsilon}(t,x))dt\nonumber\\
 &+ &\sqrt{\varepsilon}\sum_{j=1}^{\infty} h_j(t,x, u^{\varepsilon}(t,x),\nabla u^{\varepsilon}(t,x))dB^j_t\nonumber\\
 &+& \sum_{j=1}^{\infty} h_j(t,x, u^{\varepsilon}(t,x),\nabla u^{\varepsilon}(t,x))k^{\varepsilon, j}_tdt =-\nu^{\varepsilon}(dt,dx),\\
 u^{\varepsilon}(t,x)&\geq & L(t,x), \quad \quad (t,x)\in \R^+\times \R^d,\nonumber\\
 u^{\varepsilon}(T,x)&=&\Phi (x), \quad\quad x\in \R^d.
 \end{eqnarray}
 Moreover,  $v^{\varepsilon}=\Gamma^0\left(\int_0^{\cdot}k^\varepsilon(s)ds\right)$ is the solution of the random obstacle problem:
\begin{eqnarray}\label{3.2}
dv^{\varepsilon}(t,x)&+&\frac{1}{2}\Delta v^{\varepsilon}(t,x)+\sum_{i=1}^d\partial_ig_i(t,x,v^{\varepsilon}(t,x),\nabla v^{\varepsilon}(t,x))dt +f(t,x, v^{\varepsilon}(t,x),\nabla v^{\varepsilon}(t,x))dt\nonumber\\
 &+& \sum_{j=1}^{\infty} h_j(t,x, v^{\varepsilon}(t,x),\nabla v^{\varepsilon}(t,x))k^{\varepsilon, j}_tdt =-\mu^{\varepsilon}(dt,dx),\\
 v^{\varepsilon}(t,x)&\geq & L(t,x), \quad \quad (t,x)\in \R^+\times \R^d,\nonumber\\
 v^{\varepsilon}(T,x)&=&\Phi (x), \quad\quad x\in \R^d.
 \end{eqnarray}
The condition (ii) in Theorem 3.2 will be satisfied if we prove
 \begin{equation}\label{3.3}
 \lim_{\varepsilon\rightarrow 0}\left \{ E[\underset{0\leq t\leq T} {\Sup}|u^{\varepsilon}_t-v^{\varepsilon}_t|^2]+E[\int_0^T\|u^{\varepsilon}_t-v^{\varepsilon}_t\|^2dt]\right \}=0,
 \end{equation}
 here $u^{\varepsilon}_t=u^{\varepsilon}(t, \cdot)$ and $v^{\varepsilon}_t=v^{\varepsilon}(t, \cdot)$.
 The rest of the proof is to establish (\ref{3.3}).
 By Ito formula, we have
 \begin{align}\label{3.4}
\nonumber &|u^{\varepsilon}_t-v^{\varepsilon}_t|^2+\int_t^T|\nabla (u^{\varepsilon}_s-v^{\varepsilon}_s)|^2ds\\
\nonumber &=-2\int_t^T<\nabla (u^{\varepsilon}_s-v^{\varepsilon}_s), g(s,\cdot,u^{\varepsilon}_s,\nabla u^{\varepsilon}_s)-g(s,\cdot,v^{\varepsilon}_s,\nabla v^{\varepsilon}_s)>ds\\
\nonumber &+2\int_t^T<u^{\varepsilon}_s-v^{\varepsilon}_s, f(s,\cdot,u^{\varepsilon}_s,\nabla u^{\varepsilon}_s)-f(s,\cdot,v^{\varepsilon}_s,\nabla v^{\varepsilon}_s)>ds\\
\nonumber &+2\int_t^T\sum_{j=1}^{\infty}<u^{\varepsilon}_s-v^{\varepsilon}_s, h_j(s,\cdot, u^{\varepsilon}_s,\nabla u^{\varepsilon}_s)-h_j(s,\cdot, v^{\varepsilon}_s,\nabla v^{\varepsilon}_s)>k^{\varepsilon, j}_sds\\
\nonumber&+2\sqrt{\varepsilon}\sum_{j=1}^{\infty} \int_t^T<u^{\varepsilon}_s-v^{\varepsilon}_s, h_j(s,\cdot, u^{\varepsilon}_s,\nabla u^{\varepsilon}_s)>dB^j_s\\
\nonumber& +2\int_t^T<u^{\varepsilon}_s-v^{\varepsilon}_s, d\nu_s^{\varepsilon}- d\mu_s^{\varepsilon}>+\varepsilon \sum_{j=1}^{\infty} \int_t^T|h_j(s,\cdot, u^{\varepsilon}_s,\nabla u^{\varepsilon}_s)|^2ds\\
&:=I_1(t)+I_2(t)+I_3(t)+I_4(t)+I_5(t)+I_6(t).
\end{align}
Here 
$$\int_t^T<u^{\varepsilon}_s-v^{\varepsilon}_s, d\nu_s^{\varepsilon}- d\mu_s^{\varepsilon}>=\int_t^T\int_{\R^d}(u^{\varepsilon}(s,x)-v^{\varepsilon}(s,x))[\nu^{\varepsilon}(ds,dx)- \mu^{\varepsilon}(ds,dx)].$$
With the assumptions on $g$ in mind, applying Young's inequality  we have for any $\delta_1>0$
\begin{equation}\label{3.5}
I_1(t)\leq (\delta_1+2\alpha )\int_t^T|\nabla (u^{\varepsilon}_s-v^{\varepsilon}_s)|^2ds+C\int_0^t|u^{\varepsilon}_s-v^{\varepsilon}_s|^2ds.
\end{equation}
By the assumption on $f$, for any $\delta_2>0$, we have
\begin{equation}\label{3.6}
I_2(t)\leq \delta_2\int_t^T|\nabla (u^{\varepsilon}_s-v^{\varepsilon}_s)|^2ds+C\int_0^t|u^{\varepsilon}_s-v^{\varepsilon}_s|^2ds.
\end{equation}
Using the assumption on $h$, given any $\delta_3>0$, we also have
\begin{equation}\label{3.7}
I_3(t)\leq \delta_3\int_t^T|\nabla (u^{\varepsilon}_s-v^{\varepsilon}_s)|^2ds+C\int_0^t|u^{\varepsilon}_s-v^{\varepsilon}_s|^2(1+\|k^{\varepsilon}_s\|_{l^2}^2) ds.
\end{equation}
For the term $I_5$ in (\ref{3.4}), we have
\begin{equation}\label{3.8}
I_5(t)=2\int_t^T<u^{\varepsilon}_s-L(s,\cdot) + L(s,\cdot)-v^{\varepsilon}_s, d\nu_s^{\varepsilon}- d\mu_s^{\varepsilon}>\leq 0.
\end{equation}
Substituting (\ref{3.5})--(\ref{3.8}) back into (\ref{3.4}), choosing $\delta_1, \delta_2, \delta_3$ sufficiently small and rearranging terms we can find a positive constant $\delta>0$ such that
\begin{eqnarray}\label{3.9}
&&|u^{\varepsilon}_t-v^{\varepsilon}_t|^2+\delta \int_t^T|\nabla (u^{\varepsilon}_s-v^{\varepsilon}_s)|^2ds\nonumber\\
&\leq& C\int_t^T|u^{\varepsilon}_s-v^{\varepsilon}_s|^2(1+\|k^{\varepsilon}_s\|_{l^2}^2) ds\nonumber\\
&&+2\sqrt{\varepsilon}\sum_{j=1}^{\infty} \int_t^T<u^{\varepsilon}_s-v^{\varepsilon}_s, h_j(s,x, u^{\varepsilon}_s,\nabla u^{\varepsilon}_s)>dB^j_s\nonumber\\
&& +\varepsilon \sum_{j=1}^{\infty} \int_t^T|h_j(s,x, u^{\varepsilon}_s,\nabla u^{\varepsilon}_s)|^2ds.
\end{eqnarray}
By the Gronwall's inequality it follows that
\begin{eqnarray}\label{3.10}
&&\underset{0\leq t\leq T} {\Sup}|u^{\varepsilon}_t-v^{\varepsilon}_t|^2+\int_0^T\|u^{\varepsilon}_t-v^{\varepsilon}_t\|^2dt \nonumber\\
&\leq& (M_1^\varepsilon +M_2^\varepsilon) \exp(C \int_0^T(1+\|k^{\varepsilon}_s\|_{l^2}^2)ds)\leq C_M (M_1^\varepsilon +M_2^\varepsilon),
\end{eqnarray}
where
$$ M_1^\varepsilon=\underset{0\leq t\leq T}{\Sup}\left |2\sqrt{\varepsilon}\sum_{j=1}^{\infty} \int_t^T<u^{\varepsilon}_s-v^{\varepsilon}_s, h_j(s,\cdot, u^{\varepsilon}_s,\nabla u^{\varepsilon}_s)>dB^j_s\right |,$$
$$ M_2^\varepsilon= \varepsilon \sum_{j=1}^{\infty} \int_0^T|h_j(s,\cdot, u^{\varepsilon}_s,\nabla u^{\varepsilon}_s)|^2ds.$$
Using Burkholder's inequality and the boundedness of $h$, we see that
\begin{eqnarray}\label{3.11}
E[M_1^\varepsilon] &\leq &C \sqrt{\varepsilon}E[(\sum_{j=1}^{\infty} \int_0^T<u^{\varepsilon}_s-v^{\varepsilon}_s, h_j(s,\cdot, u^{\varepsilon}_s,\nabla u^{\varepsilon}_s)>^2ds)^{\frac{1}{2}}]\nonumber\\
&\leq& C \sqrt{\varepsilon} E[(\int_0^T|u^{\varepsilon}_t-v^{\varepsilon}_t|^2dt)^{\frac{1}{2}}]\nonumber\\
&&\rightarrow 0, \quad \mbox{as}\quad\quad \varepsilon\rightarrow 0,
\end{eqnarray}
where we have used the fact that $\sup_{\varepsilon}\{E[|u^{\varepsilon}_t|^2]+E[|v^{\varepsilon}_t|^2]\}<\infty$.
By the condition on $h$ in the Assumption 2.1, it is also clear that
\begin{eqnarray}\label{3.12}
E[M_2^\varepsilon] &\leq & C\varepsilon E[\int_0^T(1+ |u^{\varepsilon}_s|^2+ \| u^{\varepsilon}_s)\|^2)ds]\nonumber\\
&&\rightarrow 0, \quad \mbox{as}\quad\quad \varepsilon\rightarrow 0.
\end{eqnarray}
Assertion (\ref{3.3}) follows from (\ref{3.10})-(\ref{3.12}).

\end{document}